\documentclass[11pt]{amsart}

\usepackage{amssymb,amsthm}
\usepackage{amsmath}

\newtheorem{thm}{Theorem}
\newtheorem{theorem}{Theorem}[section]
\newtheorem{corollary}[theorem]{Corollary}

\newtheorem{lemma}[theorem]{Lemma}

\def\card#1{\vert #1 \vert}
\def\gpindex#1#2{\card {#1\colon #2}}
\def\irr#1{{\rm  Irr}(#1)}

\def\cd#1{{\rm  cd}(#1)}
\def\cent#1#2{{\bf C}_{#1}(#2)}

\begin{document}

\title[Bounding an index]{Bounding an index by the largest character degree of a $p$-solvable group}

\author{Mark L. Lewis}

\address{Department of Mathematical Sciences, Kent State University, Kent, OH 44242}

\email{lewis@math.kent.edu}


\begin{abstract}
In this paper, we show that if $p$ is a prime and $G$ is a $p$-solvable group, then $\gpindex G{O_p (G)}_p \le \left( b(G)^p/p \right)^{1/(p-1)}$ where $b(G)$ is the largest character degree of $G$.  If $p$ is an odd prime that is not a Mersenne prime or if the nilpotence class of a Sylow $p$-subgroup of $G$ is at most $p$, then $\gpindex G{O_p (G)}_p \le b(G)$.
\end{abstract}

\maketitle





\section{Introduction}

Throughout this paper, $G$ is a finite group.  Given $G$, we have the sets $\irr G$ which is the set of irreducible characters of $G$ and $\cd G = \{ \chi (1) \mid \chi \in \irr G \}$ which is the set of character degrees of $G$.  We set $b (G) = {\rm max} (\cd G)$.  It is known for solvable groups that $b (G)$ is connected with the structure of $G$.  Gluck showed in \cite{gluck} that $|G:F (G)|$ is bounded by a power of $b (G)$ where $F (G)$ is the Fitting subgroup of $G$, and Gluck conjectured when $G$ is solvable that in fact $|G:F (G)|$ is bounded by $(b (G))^2$. Gluck's conjecture has been verified for groups of odd order (see \cite{espuelas}), solvable groups whose orders are not divisible by $3$ (see \cite{yang}), and solvable groups with abelian Sylow $2$-subgroups (see \cite{doja}).  The best general result in the direction of Gluck's conjecture is that $|G:F(G)| \le (b (G))^3$ (see \cite{mowo}).  Gluck's conjecture is known to not be true for nonsolvable groups.

When we focus on a single prime, a stronger bound than the one conjectured by Gluck can be found.  This is motivated by the paper \cite{jaf}.  In \cite{jaf}, Jafari proved the following: let $G$ be a solvable group, and let $p$ be an odd prime.  If $G$ does not have a section isomorphic to $Z_p \wr Z_p$, then $\gpindex G{O_p (G)}_p \le b (G)$.  This result was an improvement on Theorem 12.29 of \cite{text} which stated that if $b (G) < p$, then $\gpindex G{O_p (G)}_p < p$ and Theorem 12.32 of \cite{text} which stated that if $b (G) < p^{3/2}$, then $\gpindex G{O_p (G)}_p < p^2$.  It also improved the result that Benjamin proved in \cite{Ben}, that if $G$ is
$p$-solvable and $b (G) < p^2$, then $\gpindex G{O_p (G)}_p < p^2$.  Benjamin also proved in her paper that if $G$ is a solvable group and $b (G) < p^\alpha$, then $\gpindex G{O_p (G)}_p < p^{2 \alpha}$, and if $\card G$ is odd, then $\gpindex G{O_p (G)}_p \le p^\alpha$.

Furthermore, in \cite{qian}, Qian showed that if $G$ is any group and $b (G) < p^2$, then $\gpindex G{O_p (G)}_p < p^2$.  In \cite{QiSh}, Qian and Shi showed that if $G$ is any group, then $\gpindex G{O_p (G)}_p < b(G)^2$ and $\gpindex G{O_p (G)}_p \le b(G)$ if $G$ is has an abelian $p$-Sylow
subgroup.

In this note, we will show that Jafari's result can be extended to $p$-solvable groups.

\begin{thm} \label{first}
Let $G$ be a $p$-solvable group, and let $p$ be a prime.  If $b (G) < p^p$, then $\gpindex G{O_p (G)}_p < p^p$.
\end{thm}

In addition, we will use Jafari's arguments to prove the following result which is along the lines of the result proved by Qian and Shi when the group has an abelian Sylow $p$-subgroup and $G$ is $p$-solvable.

\begin{thm} \label{second}
Let $G$ be a $p$-solvable group, let $p$ be a prime, and suppose that a Sylow $p$-subgroup of $G$ has nilpotence class less than $p$.  Then $\gpindex G{O_p (G)}_p \le b(G)$.
\end{thm}

Further, Jafari's arguments can also be used to prove a result along the lines of the one that Benjamin proved for groups of odd order.

\begin{thm} \label{third}
Let $G$ be a $p$-solvable group, and let $p$ be an odd prime that is not a Mersenne prime.  Then $\gpindex G{O_p (G)}_p \le b(G)$.
\end{thm}

It is easy not difficult to find examples where the bounds in Theorems \ref{first}, \ref{second}, and \ref{third} are met.  In particular, this next result shows that those bounds are optimal.

\begin{thm} \label{third examp}
Let $p$ be a prime and $n$ a positive integer.  Then there exists a solvable group $G$ with $b (G) = p^n = \gpindex G{O_p (G)}_p$ and an abelian Sylow $p$-subgroup.
\end{thm}

In general, when $p$ is a Mersenne prime and when $p = 2$, it is not true that $\gpindex G{O_p (G)}_p \le b (G)$.  We will present examples of solvable groups where $\gpindex G{O_p (G)}_p$ exceeds $b (G)$.  In our examples, $b (G) = p^{p+1}$, $G$ has a Sylow $p$-subgroup whose nilpotence class is $p$, and $p$ can be any Mersenne prime.  

\begin{thm} \label{fifth}
Let $p$ be a Mersenne prime.  Then there exists a solvable group $G$ with $b (G) = p^p$ and $\gpindex G{O_p (G)}_p = p^{p+1}$.  In particular, $\gpindex G{O_p (G)}_p = \left( b(G)^p/p \right)^{1/(p-1)}$.
\end{thm}

When $p = 2$, we do have examples with $b (G) < \gpindex G{O_2 (G)}_2$, and so, Theorem \ref{third} cannot be extended to include $p = 2$.

\begin{thm} \label{sixth}
Let $p$ be a Fermat prime so that $p = 2^f + 1$ for some positive integer $f$.  Then there exists a solvable group $G$ with $b (G) = 2^{2f}$ and $\gpindex G{O_2 (G)}_2 = 2^{2f + 1}$.
\end{thm}

However, when $p$ is either $2$ or a Mersenne prime, we are able to obtain a general bound using Jafari's argument.  To do this, we need a Theorem of Isaacs regarding the existence of a ``large'' orbit.  The result we obtain is the following:

\begin{thm} \label{fourth}
Let $G$ be a $p$-solvable group, and let $p$ be either $2$ or a Mersenne prime. Then $\gpindex G{O_p (G)}_p \le \left( b(G)^p/p \right)^{1/(p-1)}$.
\end{thm}

Notice that $\left( b(G)^p/p \right)^{1/(p-1)} < b(G)^{p/(p-1)} \le b(G)^2$, so this bound is better than the bound found by Qian and Shi for general groups and by Benjamin for $p$-solvable groups.  Of course, we are assuming the stronger hypothesis that $G$ is $p$-solvable in our result.

The examples in Theorem \ref{fifth} show that when $p$ is odd, the bound in Theorem \ref{fourth} is obtained, and so, that bound is best possible.  When $p = 2$, we do not have examples where the bound in Theorem \ref{fourth} is met, so it may be possible to improve the bound when $p = 2$.

Before we conclude this introduction, we should mention a related problem.  Define $c_p (G) = {\rm max} \{ a_p \mid a \in \cd G \}$.  In \cite{mowo}, Moret\'o and Wolf show that $\gpindex G{F(G)}_p$ is bounded in terms of $c_p (G)$.  In particular, they showed that $\gpindex G{F(G)}_p \le (c_p (G))^{19}$ when $G$ is any solvable group and that $\gpindex G{F(G)}_p \le (c_p (G))^2$ when $|G|$ is odd.  An example in \cite{esna} shows that the bound is best possible when $|G|$ is odd.  In fact, no solvable group $G$ is known to have $\gpindex G{F(G)}_p > (c_p (G))^2$, so it seems likely that the bound can be improved in general.  One might hope that the techniques used in that might apply to that problem, but at this time we have not been able to make our inductive arguments work in that setting.

We first prove our results in the solvable case.  We then present the examples.  In the final section, we prove that our bounds still hold when $G$ is a $p$-solvable group.  At this time, we have not determined whether these same bounds can be proved for any group $G$ or whether the hypothesis that $G$ is $p$-solvable is necessary to obtain our improved bounds.


\medskip

\section{Solvable Groups}

The following is essentially Theorems 4.4 and 4.8 of \cite{MaWo}.

\begin{lemma}\label{regular}
Let $P$ be a $p$-group, and assume that $P$ acts faithfully and coprimely on an abelian group $V$.  Assume one of the following conditions:
\begin{enumerate}
\item $p$ is odd and not a Mersenne prime.
\item $Z_p \wr Z_p$ is not a section of $P$.
\end{enumerate}
Then $P$ has a regular orbit on $V$.
\end{lemma}

\begin{proof}
We work by induction on $\card V$.  We may assume that $P > 1$, and since $P$ acts faithfully, this implies that $\card V > 1$.  Since the action of $P$ is faithful and coprime on $V$, we see that the action of $P$ on $V/\Phi (V)$ is faithful and coprime.  Using the inductive hypothesis, the result follows if $\Phi (V) > 1$.  Thus, we may assume that $\Phi (V) = 1$.  Therefore, we can $V$ as a completely reducible $P$-module, possibly of mixed characteristic.

If $V$ is not irreducible under the action of $P$, then we have $V = U \oplus W$ where $U$ and $W$ are nontrivial $P$-submodules of $V$. Applying the inductive hypothesis, we obtain elements $u \in U$ and $w \in W$ so that $\cent Pu = \cent PU$ and $\cent Pw = \cent PW$.  It follows that $\cent P{u+w} = \cent Pu \cap \cent Pw = \cent PU \cap \cent PW
= 1$, and we are done.  Thus, we may assume that $V$ is irreducible under the action of $P$.  If $p$ is odd and not a Mersenne prime, we can apply Theorem 4.4 of \cite{MaWo} to obtain the result.  If $Z_p \wr Z_p$ is not a section of $P$, then we can apply Theorem 4.8 of \cite{MaWo} to obtain the result.  In either case, the lemma is proved.
\end{proof}

Notice that the proof of this next theorem is essentially the proof of Theorem 1 in \cite{jaf}.  Notice that hypothesis (3) gives Theorem \ref{first} when $G$ is solvable and hypothesis (1) gives Theorem \ref{third} under the hypothesis that $G$ is solvable.

\begin{theorem}\label{Merandwre}
Let $G$ be a solvable group, and let $p$ be a prime. Assume one of the following conditions:
\begin{enumerate}
\item $p$ is odd and not a Mersenne prime.
\item $Z_p \wr Z_p$ is not a section of $G$.
\item $b (G) < p^p$.
\end{enumerate}
Then $\gpindex G{O_p (G)}_p \le b (G)$.
\end{theorem}

\begin{proof}
Assume the theorem is false and suppose that $G$ is a counterexample of minimal order.  Observe that if $G$ is a counterexample, then $G/O_p (G)$ is a counterexample.  By the minimality of $\card G$, we deduce that $O_p (G) = 1$.  Let $P$ be a Sylow $p$-subgroup of $G$, and observe that $\card P = \card G_p = \gpindex G{O_p (G)}_p$.  Let $F$ be the Fitting subgroup of $G$.  Since $G$ is a nontrivial solvable group, it follows that $F$ must be a nontrivial $p'$-group.  We know that $\cent GF \le F$.  Observe that $b (PF) \le b (G)$.  If
$FP < G$, then since $G$ is a counterexample with $\card G$ minimal, we have that $\card P < b (PF) \le  b(G)$ which is a contradiction.  Thus, we may assume that $G = FP$.  Since $P$ acts coprimely and faithfully on the nilpotent group $F$, it follows that $P$ acts coprimely and faithfully on the abelian group $F/\Phi (F)$.  Using Brauer's permutation lemma (Theorem 6.32 of \cite{text}), $P$ acts faithfully and coprimely on the abelian group $\irr {F/\Phi (F)}$.  Conditions (1) or (2) hold, then we can apply Lemma \ref{regular}.  Thus, we assume that condition (3) holds.  We can find a subgroup $Q$ in $P$ so that $\gpindex PQ = p$.  We know that $b (FQ) \le b (G) < p^p$, so we may apply the inductive hypothesis to $FQ$, to obtain the conclusion that $\card Q \le b(FQ) < p^p$.  This implies
that $\card P < p^{p+1}$, and so, $Z_p \wr Z_p$ cannot be a section of $P$ (since $\card {Z_p \wr Z_p} = p^{p+1}$).  Thus, we may apply Lemma \ref{regular} in all cases.  By Lemma \ref{regular}, there is a linear character $\lambda \in \irr {F/\Phi (F)}$, so that the stabilizer of $\lambda$ in $P$ is trivial, and thus, the stabilizer of $\lambda$ in $G$ is $F$. This implies that $\lambda^G$ is irreducible, and so, $\card P = \gpindex GF = \lambda^G (1) \le b(G)$.  This proves the theorem.
\end{proof}

As an immediate corollary, we obtain Theorem \ref{second} under the solvable hypothesis.

\begin{corollary}
Let $G$ be a solvable group, and let $p$ be a prime.  Assume that a Sylow $p$-subgroup of $G$ has nilpotence class less than $p$.  Then $\gpindex G{O_p (G)}_p \le b (G)$.
\end{corollary}

\begin{proof}
Notice that if $Z_p \wr Z_p$ is a section of $G$, then it must be a section of some Sylow $p$-subgroup of $G$.  This would imply that a Sylow $p$-subgroup of $G$ has nilpotence class at least $p$ which is a contradiction.  Thus, $Z_p \wr Z_p$ is not a section of $G$, and we can apply Theorem \ref{Merandwre} to obtain the result.
\end{proof}

Using the following result of Isaacs, we can remove the extra
hypotheses on $p$ or on the structure of a Sylow subgroup used in the previous results.  Unfortunately, we have to weaken our
conclusion in this case, and we will present examples to show that the weaker conclusion is necessary.  The following is proved as Theorem A in \cite{large}.  While it does not prove that there is a regular orbit, it does prove that there is a ``large'' orbit.

\begin{theorem}\cite{large}\label{large}
Let $P$ be a $p$-group that acts faithfully and coprimely on a group $V$.  Then there exists an element $v \in V$ so that $\card {\cent Pv} \le \left( \card P/p \right)^{1/p}$.
\end{theorem}

With this result in hand, we can prove the following which includes Theorem \ref{fourth} when $G$ is solvable.

\begin{theorem}\label{general}
Let $G$ be a solvable group and let $p$ be a prime.  Then $\gpindex G{O_p (G)}_p \le \left( b (G)^p/p \right)^{1/(p-1)}$.
\end{theorem}

\begin{proof}
We work by induction on $\card G$.  Using the inductive hypothesis, we may assume that $O_p (G) = 1$. Let $P$ be a Sylow $p$-subgroup of $G$, and note that $\card P = \card G_p = \gpindex G{O_p (G)}_p$.  Let $F$ be the Fitting subgroup of $G$. Since $G$ is solvable, we know that $\cent GF \le F$ and $F$ is a $p'$-group. Observe that $b (PF) \le b (G)$. If $PF < G$, then we obtain the result by the inductive hypothesis. Hence, we may assume $G = PF$.  Now, $P$ acts faithfully and coprimely on the nilpotent group $F$.  Thus, $P$ will
act faithfully and coprimely on $F/\Phi (F)$. Applying Theorem
\ref{large}, we can find a linear character $\lambda \in \irr
{F/\Phi (F)}$ so that $\card {\cent P{\lambda}} \le \left( \card P/p \right)^{1/p}$. Thus, if $T$ is the stabilizer of $\lambda$ in $G$, then $T = F \cent P{\lambda}$, and so
$$
\gpindex GT = \gpindex P{\cent P{\lambda}} \ge \frac {\card
P}{\left( \frac {\card P}p \right)^{1/p}} = \left( \card P^{p-1} p \right)^{1/p}.
$$
By Clifford's theorem, it follows that $\gpindex GT \le b (G)$.  This implies that $\left( \card P^{p-1} p \right)^{1/p} \le b (G)$, and we conclude that $\card P \le \left( b(G)^p/p \right)^{1/(p-1)}$, as desired.
\end{proof}

\section{Examples}

In this section, we present three examples.  In particular, we prove Theorems \ref{third examp}, \ref{fifth}, and \ref{sixth}. We first present the proof of Theorem \ref{third examp}

\begin{proof}[Proof of Theorem \ref{third examp}]
Let $p$ be a prime and let $n$ be a positive integer.  Let $q$ be a prime that is different from $p$.  In particular, $q$ is relatively prime to $p$, and so, $q$ is a unit modulo $p$.  It follows that there is a positive integer $m$ so that $p^n$ divides $q^m - 1$.  Let $F$ be the Galois field of order $q^m$ and let $C$ be the subgroup of order $p^n$ in the multiplicative group of $F$.  It is clear that $C$ acts via automorphisms on the additive group of $F$, and take $G$ to be the resulting semi-direct product.  It is easy to see that $b (G) = p^n$ and $\gpindex G{O_p (G)}_p = p^n$.  
\end{proof}

We now present examples to see that the additional hypothesis is needed in Theorem \ref{Merandwre} and to see that the bound in Theorem \ref{general} is appropriate.  The next two examples can be found as Example 4.5 of \cite{MaWo}.  Thus, the bound in Theorem \ref{general} is optimal.

\begin{proof}[Proof of Theorem \ref{fifth}]
Let $p$ be a Mersenne prime, so $p = 2^f - 1 \ge 3$.  We take $P = Z_p \wr Z_p$, and we take $V$ to be $p$ copies of ${\rm GF} (2^f)$. In \cite{MaWo}, they show that $P$ acts faithfully on $V$, and that $P$ does not have any regular orbits in its action on $V$.  Let $G$ be the semi-direct product of $P$ acting on $V$. We claim that $b (G) = p^p$. On the other hand, $O_p (G) = 1$, so $\gpindex G{O_p (G)}_p = \card G_p = \card P = p^{p+1}$ which is larger than $b (G)$.  In addition, $\left( b(G)^p/p \right)^{1/(p-1)} = \left ( (p^p)^p/p \right)^{1/(p-1)} = (p^{p^2-1})^{1/(p-1)} = p^{p+1}$.   This proves Theorem \ref{fifth}.
\end{proof}

Finally, we have the proof of Theorem \ref{sixth}.

\begin{proof}[Proof of Theorem \ref{sixth}]
In this example, we take $p = 2$.  Let $q$ be a Fermat prime, so $q = 2^f + 1 \ge 3$.  Let $P = Z_{2^f} \wr Z_2$ and let $V$ be $2$ copies of ${\rm GF} (q)$.  In \cite{MaWo}, they show that $P$ acts faithfully on $V$, and has no regular orbit. Let $G$ be the semi-direct product of $P$ acting on $V$. We claim that $b (G) = 2^{2f}$.  On the other hand, $O_2 (G) = 1$, so $\gpindex G{O_2 (G)}_2 = \card G_2 = \card P = 2^{2f+1}$.  In particular, $\gpindex G{O_2 (G)}_2$ is large than $b (G)$.  This yields Theorem \ref{sixth}.
\end{proof}

\section{$p$-solvable groups}

This next lemma follows from a corollary of Gluck's regular orbit theorem.

\begin{lemma}[Corollary 5.7(b) of \cite{MaWo}] \label{dolfi}
Let $p$ be a prime an odd prime, and suppose that $P$ is a $p$-group that is a permutation group on $\Omega$.  Then there exists a set $\Delta \subseteq \Omega$ so that $P_{\Delta} = 1$.
\end{lemma}

We also need the following lemma which is proved in \cite{MoTi}.

\begin{lemma}[Proposition 2.6 of \cite{MoTi}] \label{MoTi}
Let $A$ act faithfully and coprimely on a nonabelian simple group $S$.  Then $A$ has at least $2$ regular orbits on $\irr S$.
\end{lemma}

Using these two results, we obtain a character that induces
irreducibly in the key situation where we have a $p$-group that is acting on a direct product of copies of a nonabelian simple group.

\begin{lemma} \label{induction}
Let $S$ be a nonabelian simple group, and let $p$ be a prime that does not divide $\card S$.  Suppose $V = S_1 \times \cdots \times S_n$ where $S_i \cong S$.  Assume $P$ is a $p$-group that acts faithfully on $V$ via automorphisms, and assume the action of $P$ transitively permutes the $S_i$'s.  If $G$ is the semi-direct product of $P$ acting on $V$, then there exists $\theta \in \irr V$ so that $\theta^G$ is irreducible.
\end{lemma}

\begin{proof}
Let $\Omega = \{ S_1, \dots, S_n \}$.  Let $Q$ be the kernel of the action of $P$ on $\Omega$.  Notice that since $p$ does not divide $\card S$, it follows that $p$ is odd.  Thus, we can use Lemma \ref{dolfi} to obtain a set $\Delta \subseteq \Omega$ so that $\cent P{\Delta} = Q$. Since $P$ acts transitively on $\Omega$, it follows that $P$ acts transitively on the set $\{ \cent Q{S_1}, \dots, \cent Q{S_n} \}$. Let $R_i = Q/\cent Q{S_i}$, and this implies that there is a group $R$ so that $R \cong R_i$ for all $i$.  Notice that $R$ acts on $S$.  By Lemma \ref{MoTi}, we know that $R$ has two regular orbits on $\irr S$.  Hence, we can find $\mu, \nu \in \irr S$ in different $R$-orbits so that $\cent P{\mu} = \cent P{\nu} = 1$. For each $i$, let $\mu_i$ and $\nu_i$ be the characters in $\irr S$ that correspond to $\mu$ and $\nu$ respectively.

We now define $\theta \in \irr V$ as follows.  For each $i$, we take $\theta_i = \mu_i$ if $S_i \in \Delta$ and $\theta_i = \nu_i$ if $S_i \not\in \Delta$, and then, $\theta = \theta_1 \times \cdots \times \theta_n$.  Observe that $\cent P{\theta} \le \cent P{\Delta} = Q$, and so,
$$
\cent P{\theta} = \cent Q{\theta} = \cap_{i=1}^n \cent Q{\theta_i} = \cap_{i=1}^n \cent Q{S_i} = 1.
$$
This implies that the stabilizer of $\theta$ in $G$ is $V$.  Hence, $\theta^G$ is irreducible.
\end{proof}

We now prove that the results that we proved when $G$ is solvable still hold when $G$ is $p$-solvable.  To simplify the next result we define the following function $f$.  We define $f (G,p) = b(G)$ if $p$ is odd and not a Mersenne prime or if $p= 2$ or $p$ is a Mersenne prime and $G$ has no section isomorphic to $Z_p \wr Z_p$.  We define $f (G,p) = \left( b(G)^p/p \right)^{1/(p-1)}$ if either $p = 2$ or $p$ is a Mersenne prime and $G$ has a section isomorphic to $Z_p \wr Z_p$.  Thus, this next theorem includes the $p$-solvable case of all of the theorems in the Introduction.

\begin{theorem}
Let $G$ be a $p$-solvable group where $p$ is some prime. Then
$\gpindex G{O_p (G)}_p \le f (G,p)$.
\end{theorem}

\begin{proof}
We will work by induction on $\card G$.  Suppose $H$ is a section of $G$.  We claim that $f (H,p) \le f (G,p)$.  If $p$ is odd and not a Mersenne prime, then $f (H,p) = b (H) \le b (G) = f (G,p)$.  Suppose for now that $p = 2$ or $p$ is a Mersenne prime.  If $G$ does not have a section isomorphic to $Z_p \wr Z_p$, then neither does $H$, and again, we have $f (H,p) = b (H) \le b (G) = f (G,p)$.  If $H$ has a section isomorphic to $Z_p \wr Z_p$, then so does $G$, and so,
$f (H,p) = \left( b (H)^p/p \right)^{1/(p-1)} \le \left( b (G)^p/p \right)^{1/(p-1)} = f (G,p)$.  Finally, we must consider the case where $G$ has a section isomorphic to $Z_p \wr Z_p$ and $H$ does not.  In this case, since $Z_p \wr Z_p$ is a section of $G$, it follows that $G$ must have an irreducible character whose degree is at least $p$, and hence, $p \le b (G)$.  Thus, $b (G)^{(p-1)} p \le b (G)^p$, and $b (G) \le \left( b (G)^p/p \right)^{1/(p-1)}$.  We deduce that $f (H,p) = b (H) \le b (G) \le \left( b (G)^p/p \right)^{1/(p-1)} = f (G,p)$.  This proves the claim in all cases.

Using the inductive hypothesis on $G/O_p (G)$, we may assume that $O_p (G) = 1$. Also, using the Frattini argument, one can now show that $O_p (G/\Phi (G)) = 1$, and so, using the inductive hypothesis on $G/\Phi (G)$, we may assume that $\Phi (G) = 1$.  If all the minimal normal subgroups of $G$ are solvable, then let $F$ be the Fitting subgroup of $G$, and it is known that there is a subgroup $A$ so that $G = FA$ and $F \cap A = 1$.  It is easy to see that $\cent AF$ will be a normal subgroup of $G$, and since $F$ contains all the minimal normal subgroups of $G$, we conclude that $\cent AF = 1$, and hence, $\cent GF \le F$.  So it suffices to show that the
result holds in $FP$, but $FP$ is solvable, and we have seen that the result holds in solvable groups.

Thus, we may assume that $G$ has a nonsolvable minimal normal
subgroup $V$.  Since $G$ is $p$-solvable, we see that $p$ does not divide $\card V$.  Notice that this implies that $p$ must be odd. We know that $\cent GV$ is normal in $G$, so $O_p (\cent GV) \le O_p (G) = 1$. Let $H = V \cent GV P$. Now, $O_p (H)$ and $V$ are normal in $H$ and of coprime orders, so $O_p (H) \le \cent GV$, and this implies that $O_p (H) = 1$. Thus, we may use the inductive hypothesis to assume that $G = H = V \cent GV P$. Observe that $V = S_1 \times \cdots \times S_n$ where $S_i \cong S$ and $S$ is a nonabelian simple group.  This implies that $V \cap \cent GV = 1$.  

By Lemma \ref{induction}, we can find $\theta \in \irr V$ so that $\cent P{\theta} = \cent PV$. Let $\gamma \in \irr {\cent GV}$ so that $\gamma (1) = b (\cent PV)$. Notice that the stabilizer of $\theta \times \gamma$ is $V \cent GV$.  This implies that $(\theta \times \gamma)^G$ is irreducible, and hence, $b (G) > \gpindex
P{\cent PV} b (\cent GV)$. By the inductive hypothesis applied in $\cent GV$, we have that $\card {\cent PV} \le f (\cent GV, p)$. If either $p$ is not a Mersenne prime or $Z_p \wr Z_p$ is not involved in $G$, then $f (G,p) = b (G)$ and $f (\cent GV, p) = b (\cent GV)$, and the result follows.  Thus, we may assume that $p$ is a Mersenne prime and $Z_p \wr Z_p$ is involved in $G$.  Then we have that
$$
f (G,p) = \left( b(G)^p/p \right)^{1/(p-1)} > \left ( \gpindex
P{\cent PV}^p b(\cent GV)^p/p \right)^{1/(p-1)}.
$$
This expression can be rewritten as follows: 
$$
\gpindex P{\cent PV}^{p/(p-1)} \left( b(\cent GV)^p/p \right)^{1/(p-1)}.
$$
Using the definition of $f (\cent GV,p)$, and the fact that $p/(p-1) > 1$, we see that this is larger than $\gpindex P{\cent PV} f (\cent GV,p)$, and obtain
$$
f (G,p) > \gpindex P{\cent PV} f (\cent GV,p) \ge \gpindex P{\cent PV} \card {\cent PV} = \card P.
$$
This proves the theorem.
\end{proof}

\end{document}